\documentclass[11pt]{article}
\usepackage{latexsym, amsmath, amssymb} 

\usepackage{amsmath}
\usepackage{amsfonts}
\usepackage{amssymb}
\usepackage{amscd}
\usepackage{amsthm}
\usepackage{fancyhdr}
\theoremstyle{plain}
\newtheorem{theorem}{Theorem}[section]
\newtheorem{proposition}[theorem]{Proposition}

\newtheorem{lemma}[theorem]{Lemma}
\newtheorem{definition}[theorem]{Definition}
\newtheorem{corollary}[theorem]{Corollary}

\newcommand\sD{{\mathcal D}}

\newcommand\sF{{\mathcal F}}
\newcommand\sG{{\mathcal G}}
\newcommand\sH{{\mathcal H}}
\newcommand\sI{{\mathcal I}}

\newcommand\sO{{\mathcal O}}

  \def \tab#1{\kern #1 truein}

\begin{document}
 \title{THE EQUATIONS OF SINGULAR LOCI \\OF AMPLE DIVISORS ON \\(SUBVARIETIES OF) ABELIAN VARIETIES}
\author{Luigi Lombardi and Francesco Malaspina}
\maketitle

  \begin{abstract}
In this paper we consider ideal sheaves associated to  the singular loci
   of a divisor in a linear system $|L|$  of an ample line bundle on a complex abelian variety. We prove an effective result on their (continuous) global generation, after suitable twists by powers of $L$. 
    Moreover we show that similar results hold for subvarieties of a complex
abelian variety.
\end{abstract}

\section{Introduction}

In the paper \cite{PP2} Proposition $7.21$, G. Pareschi and M. Popa studied the equations of the special subvarieties $W_d$ in Jacobians by means of  theta-regularity and the continuous global generation of sheaves on abelian varieties. In the same vein, they gave an effective bound for the equations of the singular locus of the theta divisor on  a Jacobian, $\Sigma (\Theta) \cong W^1_{g-1}$, by showing that the ideal sheaf $\sI_{\Sigma (\Theta)}$ is $3$-$\Theta$-regular (cf. \cite{PP2}, Proposition $7.21$) and hence $\sI_{\Sigma (\Theta)} (3\Theta)$ is globally generated.\\ In this paper we generalize this result to an arbitrary complex abelian variety, where we consider an arbitrary ample line bundle instead of a theta divisor. Moreover we consider the multiplicity-$k$ locus of a divisor $D\in |L|$. Now we present the generalization.

Let $A$ be a complex abelian variety of dimension $g\geq 2$ and $L$ be an ample line bundle on $A$. Let $D\in |L|$ be a divisor, $$\Sigma_k(D)=\{x\in A\;|\;\mbox{mult}_xD>k\},$$ be the multiplicity-$k$ locus of $D$ and $\sI_{\Sigma_k(D)}$ be its ideal sheaf:

\begin{theorem}\label{mainth}
 The sheaf $\sI_{\Sigma_k(D)}\otimes L^{\otimes 3}$ is globally generated.
\end{theorem}

This Theorem allows us to find the degrees of defining equations of $\Sigma (D)$. In particular, it implies that $\Sigma_k(D)$ is \textit{cut-out by equations in} $|L^{\otimes 3}|$, i.e. locally there exist divisors $D_1,\ldots ,D_m\in |L^{\otimes 3}|$ such that $\Sigma_k(D)=D_1\cap \ldots \cap D_m$.

A more general problem in the case of Jacobians is to find for which positive integers $k$ the sheaves $\sI_{W_d^r(C)}(k\Theta)$ are globally generated, where $W_d^r(C)$ are the Brill-Noether loci on a smooth curve $C$, i.e. $$W_d^r(C)=\{L\in \mbox{Pic}^d(C)\;|\;h^0(L)\geq r+1\}.$$ Let $\Theta$ be a theta divisor on the jacobian $J(C)$ of a smooth curve $C$ of genus $g$. Via the identification $$\Theta\cong W_{g-1}^0(C)=\{L\in \mbox{Pic}^{g-1}(C)\;|\;h^0(L)\geq 1\},$$ Riemann's Theorem ensures that $$\Sigma_k(\Theta)=W_{g-1}^k(C),$$ so Theorem \ref{mainth} is a result in this direction.

The notion of globally generated sheaf is not the unique way in order to get equations of a subvariety. Let $A$ be a complex abelian variety and $\sF$ be a sheaf on $A$. We say that the sheaf $\sF$ is \textit{continuously globally generated} if for any non-empty open subset $U\subset \mbox{Pic}^0(A)$ the sum of evaluation maps $$\bigoplus_{\alpha \in U} H^0(\sF \otimes \alpha)\otimes \alpha^\vee \longrightarrow \sF$$ is surjective (see \cite{PP1}). With the same notation as in Theorem \ref{mainth} we have another result:
\begin{theorem}\label{mainth2}
 The sheaf $\sI_{\Sigma_k(D)}\otimes L^{\otimes 2}$ is continuously globally generated.
\end{theorem}

We will see that Theorem \ref{mainth} is an easy consequence of Theorem \ref{mainth2}.
In particular, Theorem \ref{mainth2} implies that $\Sigma_k(D)$ is \textit{cut-out by equations in} $|L^{\otimes 2}\otimes \alpha|$, for some  $\alpha\in \mbox{Pic}^0(A)$, i.e. there exist line bundles $\alpha_1,\ldots ,\alpha_t\in \mbox{Pic}^0(A)$ and divisors $D_1,\ldots ,D_m\in \bigcup _i|L^{\otimes 2}\otimes \alpha_i|$ such that locally $\Sigma_k(D)=D_1\cap \ldots \cap D_m$.

The proofs of Theorem \ref{mainth} and Theorem \ref{mainth2} use a general method different from the ad-hoc argument in \cite{PP2}, Proposition $7.21$. Our main tool is the use of the bundle of differential operators associated to an ample line bundle on a complex abelian variety, see \cite{ELN}. In this case we will see that the bundle of differential operators satisfies nice cohomological properties.

In the last section we investigate the same problem on subvarieties of a complex abelian variety. More precisely, let $X\subset A$ be a complex projective smooth subvariety of dimension $n\geq 2$ of a complex abelian variety $A$, and let $M$ be an ample line bundle on $X$ and $D\in |M\otimes \omega_X|$ be a divisor. Putting $L:=M\otimes \omega_X$ we have the following results:

\begin{theorem}\label{intr-sub}
\item (i.) The sheaf $\sI_{\Sigma_1(D)}\otimes L^{\otimes 2}\otimes \omega_X$ is continuously globally generated.
\item (ii.) The sheaf $\sI_{\Sigma_1(D)}\otimes L^{\otimes 3}\otimes \omega_X$ is globally generated.
\item (iii.) The sheaf $\sI_{\Sigma_1(D)}\otimes L^{\otimes {n+2}}$ is continuously globally generated and
\item (iv.) The sheaf $\sI_{\Sigma_1(D)}\otimes L^{\otimes {n+3}}$ is globally generated.
\end{theorem}

In order to prove the last two points we will remark that the cotangent bundle of a subvariety of an abelian variety is nef and  we will state a vanishing theorem for varieties with nef cotangent bundle.

\section{Notations and Preliminaries}

Throughout this paper every variety is assumed to be irreducible. If $Y$ is a subvariety, its ideal sheaf is denoted by $\sI_Y$.

In this section we present the notion of sheaf satisfying the Index Theorem, that is a condition on the cohomology of the sheaf. After, following \cite{PP1}, we will give the definition of continuously globally generated sheaf, putting them in relation with globally generated sheaves and sheaves satisfying the index theorem. Only in this section every variety is defined over an algebraically closed field of arbitrary characteristic.

\begin{definition}[\textbf{Sheaf Satisfying the Index Theorem with Index $i$}]
A sheaf $\sF$ on an abelian variety $A$ satisfies the \emph{index theorem with index} $i$, I.T. $i$ for short, if $$H^j(\sF \otimes \alpha )=0$$ for any $\alpha\in \mbox{Pic}^0(A)$ and for any $j\neq i$.
\end{definition}

An ample line bundle on an abelian variety satisfies I.T. $0$: see for example \cite{MU} Application I, p.60, and Chapter 16. In characteristic zero it is a simple consequence of Kodaira's Vanishing Theorem.\\
Recall that a sheaf $\sF$ on an abelian variety $A$ is \emph{globally generated} if the evaluation map $H^0(\sF)\otimes \sO_A\rightarrow \sF$ is surjective. A similar notion is the following

\begin{definition}[\textbf{Continuously Globally Generated Sheaf}]
A sheaf $\sF$ on an abelian variety $A$ is \emph{continuously globally generated} if for any non-empty open subset $U\subset \mbox{Pic}^0(A)$ the sum of evaluation maps $$\bigoplus_{\alpha \in U} H^0(\sF \otimes \alpha)\otimes \alpha^\vee \longrightarrow \sF$$ is surjective.
\end{definition}



 The link between this kind of sheaves and globally generated  sheaves is explained by the following two propositions.

\begin{proposition}\label{cgg per cgg=gg}
Let $\sF$ be a coherent continuously globally generated sheaf on an abelian variety $A$ and $H$ be a continuously globally generated sheaf on $A$ which is everywhere of rank one on its support, then $\sF\otimes H$ is globally generated.
\end{proposition}
\begin{proposition}\label{it0 implica cgg}
 If $\sF$ is a sheaf satisfying I.T. $0$ on an abelian variety $A$, possibly supported on a subvariety $X$ of $A$, then $\sF$ is continuously globally generated.
\end{proposition}

For the proof of Proposition \ref{cgg per cgg=gg} see Lemma 2.3 in \cite{PP3} and for the proof of Proposition \ref{it0 implica cgg} see Proposition 2.13 in \cite{PP1} where it is stated in a more general setting.

In the sequel we will use the following Lemma.

\begin{lemma}\label{cgg-prop}
Let $A$ be an abelian variety.
\item[](i.) A quotient of a continuously globally generated sheaf $\sG$ on $A$ is still a continuously globally generated sheaf.

 \item[] (ii.) Let $$0\longrightarrow \sF'\longrightarrow \sF \longrightarrow \sF''\longrightarrow 0$$ be an exact sequence of sheaves on $A$, where $\sF'$ and $\sF''$ are continuously globally generated sheaves and such that $H^1(\sF'\otimes \alpha)=0$ for any $\alpha \in \mbox{Pic}^0(A)$. Then the sheaf $\sF$ is continuously globally generated.

\begin{proof}
 \item (i.) Let $\sG'$ be a quotient of $\sG$ and $U\subset \mbox{Pic}^0(A)$ be a non-empty open subset. Consider the following commutative diagram
$$\begin{array}{ccc}
 \bigoplus_{\alpha \in U} H^0(\sG \otimes \alpha)\otimes \alpha^\vee & \stackrel{\sigma}{\longrightarrow}& \sG
\\
 \downarrow  & &\downarrow \nu \\
\bigoplus_{\alpha \in U} H^0(\sG' \otimes \alpha)\otimes \alpha^\vee &\stackrel{\tau}{\longrightarrow}& \sG'.\\

 \end{array}$$

Since the maps $\sigma$ and $\nu$ are surjective, $\tau$ also has to be surjective, thus $\sG'$ is continuously globally generated.

\item (ii.) If $\sH$ is a sheaf and $U\subset \mbox{Pic}^0(A)$ is a non-empty subset, denote by $\bar{\sH}_U$ the sheaf $\bigoplus_{\alpha\in U} H^0(\sH\otimes \alpha)\otimes \alpha^\vee$. The hypotheses imply that for any non-empty open subset $U\subset \mbox{Pic}^0(A)$ there is a commutative diagram
$$\begin{array}{ccccccccc}
 0 & \longrightarrow & \bar{\sF'_U} & \longrightarrow & 
\bar{\sF}_U& \longrightarrow & \bar{\sF''_U} & \longrightarrow & 0\\

 & & \downarrow & & \downarrow & & \downarrow & &\\
0 & \longrightarrow & \sF'& \longrightarrow & \sF & \longrightarrow & \sF''& \longrightarrow & 0\\

 \end{array}$$
where the first and the third vertical arrow are surjective. At this point the five-lemma implies that the middle vertical arrow is surjective as well.
\end{proof}
\end{lemma}

\section{The Bundle of Differential Operators}

In order to prove Theorem \ref{mainth} and Theorem \ref{mainth2},  we need to introduce the bundle of differential operators of order $\leq k$, see \cite{ELN}, \cite{DS} and \cite{MA}.

Let $X$ be a smooth complex projective variety of dimension $n$ and $L$ be a line bundle on $X$. Let $\Delta \subset X\times X$ be the diagonal of $X$ and $p,q:X\times X\rightarrow X$ be the two projections onto the first and second factor. The \textit{k-jet bundle associated to $L$}, $J_k(L)$, is the vector bundle $$p_*(\sO_{X\times X}/\sI_{\Delta}^{k+1}\otimes q^*L)$$ where $$\sI_{\Delta}^{k+1}=\{f\in \sO_{X\times X}\; | \; \mbox{ord}_x(f)\geq k+1 \;\;\mbox{for any}\;x\in \Delta\}.$$ The k-jet bundle is a vector bundle of rank $\binom{k+n}{n}$ whose fiber is $$(J_k(L))_x=L_x\otimes \sO_{X,x}/ \frak{m}_x^{k+1},$$ where $x\in X$ and $\frak{m}_x$ is the maximal ideal of $\sO_{X,x}$. In other words the elements of a fiber are equivalence classes of sections of $L$, where two sections are in the same class if their Taylor expansions coincide up to order $k$ near $x$.

There are natural maps of sheaves $$j^k:L\longrightarrow J_k(L)$$ sending the germ of a section $s$ at a point $x\in X$ to its $k$-th jet. More specifically for $s\in H^0(X,L)$, $j^k(s(x))$ is the $\binom{k+n}{n}$-ple determined by the coefficients of the terms of degree up to $k$, in the Taylor expansion of $s$ around $x$.

In this way we get a natural projection map $$J_k(L)\longrightarrow J_{k-1}(L).$$ A germ of a section of $J_k(L)$ at a point $x\in X$ is sent to zero, under the projection map, if the terms of degree up to $k-1$ in its Taylor expansion vanish, hence the kernels of these maps are the vector bundles $\mbox{Sym}^k(\Omega_X^1)\otimes L$. In fact a germ of a section of $\mbox{Sym}^k(\Omega_X^1)\otimes L$ at a point $x\in X$ corresponds to the $\binom{k+n-1}{n-1}$-ple determined by the coefficients of the terms of degree $k$ in its Taylor expansion around $x$.

Thus, for $k\geq 1$, there are exact sequences of sheaves of $\sO_X$-modules $$0\longrightarrow \mbox{Sym}^k(\Omega ^1_X)\otimes L\longrightarrow J_k(L)\longrightarrow J_{k-1}(L)\longrightarrow 0.$$

Now we define the \textit{bundle of differential operators of order $\leq k$ associated to $L$} as $$\sD_L^k:=\mathcal{H}\textit{om}_{\sO_X}  (J_k(L),L)=J_k(L)^{\vee}\otimes L.$$ By dualizing and after tensoring by $L$ the previous exact sequences we get new exact sequences of sheaves of $\sO_X$-modules.
\begin{equation}\label{2-succ-esa}
0\longrightarrow \sD_L^{k-1}\longrightarrow \sD^k_L \longrightarrow \mbox{Sym}^k(TX)\longrightarrow 0.
\end{equation}
A non-zero section $s\in H^0(X,L)$ determines a morphism of vector bundles $$d_k(s):\sD^k_L \longrightarrow L$$ in this way. Let $U\subset X$ be an open subset and let $f:=s_{|U}\in L(U)$, the map associates to any differential operator $\Psi_U\in \sD_L^k(U)$ the section $\Psi_U(j^k_U(f))\in L(U)$, where $L(U)$ and $\sD_L^k(U)$ are the rings of sections of $L$ and of $\sD_L^k$ over the open set $U$ and $j^k_U$ is the map $j^k$ on the open set $U$.
It follows that $d_k(s)$ is zero exactly at the locus where $s$ vanishes to order $>k$. More precisely let $D\in |L|=\mathbb{P}(H^0(X,L))$, then $D$ corresponds to a section (modulo scalars), say $\phi$. Consider the multiplicity-$k$ locus $$\Sigma_k(D)=\{x\in X\; |\; \mbox{mult}_x D>k\},$$ with its natural scheme structure. Then the image of $d_k(\phi)$ is just the ideal sheaf of this scheme, i.e. one has a surjective sheaf morphism
\begin{equation}\label{sur-imp}
\sD_L^k\longrightarrow \sI_{\Sigma_k(D)}\otimes L.
\end{equation}

\section{Proof of the Theorem}

Now we are ready to prove the following

\begin{theorem}\label{mainth-proof}
 Let $A$ be a complex abelian variety of dimension $g\geq 2$ and $L$ be an ample line bundle on $A$. For any $k\geq 1$ and any divisor $D\in |L|$, let $\Sigma_k(D)=\{x\in A \; | \; \mbox{mult}_x D>k\}$ be the  multiplicity-$k$ locus of $D$. Then the sheaf $\sI_{\Sigma_k(D)}\otimes L^{\otimes s}$ is continuously globally generated for any $k\geq 1$ and any $s\geq 2$.

\begin{proof}
Let $J_k(L)$ be the $k$-jet bundle associated to $L$ and $\sD_L^k$ its bundle of differential operators. Note that $J_0(L)=L$ and that $\sD_L^0=\sO_A$.
First of all we will prove that the bundle $\sD_L^k\otimes L^{\otimes {s-1}}$ satisfies I.T.0. This is done by induction  on $k$, we begin with $k=1$.

Since on an abelian variety the tangent bundle is trivial, by the exact sequence (\ref{2-succ-esa}) we get a new exact sequence $$0\longrightarrow \sO_A\longrightarrow \sD_L^1\longrightarrow \bigoplus_g \sO_A\longrightarrow 0,$$ and by tensoring it by $L^{\otimes {s-1}}$, with $s\geq 2$, we get $$0\longrightarrow L^{\otimes {s-1}}\longrightarrow \sD_L^1\otimes L^{\otimes {s-1}}\longrightarrow \bigoplus_g L^{\otimes {s-1}}\longrightarrow 0.$$
The line bundle $L$ satisfies I.T. $0$ since it is ample on an abelian variety, for the same reason also $L^{\otimes {s-1}}$ and $\bigoplus_g L^{\otimes {s-1}}$ satisfy I.T. $0$. Therefore $\sD_L^1\otimes L^{\otimes {s-1}}$ satisfies I.T. $0$.

Now suppose that the bundle $\sD_L^{k-1}\otimes L^{\otimes {s-1}}$ satisfies I.T. $0$. By the exact sequence (\ref{2-succ-esa}) and by tensoring it by $L^{\otimes {s-1}}$ we get a new exact sequence $$0\longrightarrow \sD_L^{k-1}\otimes L^{\otimes {s-1}}\longrightarrow \sD_L^{k}\otimes L^{\otimes {s-1}}\longrightarrow \bigoplus_{\binom{g+k-1}{g-1}} L^{\otimes {s-1}}\longrightarrow 0.$$
The first term of the sequence satisfies I.T. $0$ by inductive hypothesis and easily also the last term satisfies I.T. $0$, hence also the middle term satisfies I.T. $0$. By Proposition \ref{it0 implica cgg} the bundle $\sD_L^k\otimes L^{\otimes {s-1}}$ is also continuously globally generated.

Let $D\in |L|$. By tensoring the surjection (\ref{sur-imp}) by $L^{\otimes {s-1}}$, we get a new surjection $$\sD_L^k\otimes L^{\otimes {s-1}}\longrightarrow \sI_{\Sigma_k({D})}\otimes L^{\otimes s},$$ and hence the quotient $\sI_{\Sigma_k({D})}\otimes L^{\otimes s}$ is continuously globally generated by \\Lemma \ref{cgg-prop} \emph{(i.)}.
\end{proof}
\end{theorem}
With the same notation as in the previous Theorem, we have the following

\begin{corollary}\label{last}
 The sheaf $\sI_{\Sigma_k(D)}\otimes L^{\otimes s}$ is globally generated for any $k\geq 1$ and any $s\geq 3$.

\begin{proof}
 By Theorem \ref{mainth-proof} we have that the sheaf $\sI_{\Sigma_k(D)}\otimes L^{\otimes s}$ is continuously globally generated for any $k\geq 1$ and any $s\geq 2$. By tensoring this sheaf by $L$, which is continuously globally generated by Proposition \ref{it0 implica cgg}, we get the claimed result by Proposition \ref{cgg per cgg=gg}.
\end{proof}
\end{corollary}

Finally, we would like to point out that, while the basic properties of continuously globally generated sheaves hold in arbitrary characteristic, those of the bundles of differential operators hold only in characteristic zero, for differentiation reasons. Hence the proof of Theorem \ref{mainth-proof} given above does not work in positive characteristic.

\section{Subvarieties of an Abelian Variety}

In this section we investigate the same problem on smooth subvarieties of a complex abelian variety.

\begin{definition}[\textbf{Nef Bundles}]
A line bundle $L$ on a projective variety $X$ is \textit{nef} (or \textit{numerically effective}) if for every curve $C\subset X$ $$\int_C c_1(L)\geq 0.$$

A vector bundle $E$ on a projective variety $X$ is \textit{nef} (or \textit{numerically effective)} if the associated line
 bundle $\sO_{\textbf{P}(E)}(1)$ is nef on the projectivized bundle $\textbf{P}(E)=Proj(\bigoplus _m\mbox{Sym}^mE)$.
\end{definition}

For generalities on nef vector bundles see for example \cite{LA} Theorem 6.2.12. Recall that quotients and pull-backs of
 nef vector bundles are nef, and that any tensor product, exterior product, symmetric product, direct sums and extensions of nef
 bundles are again nef. Moreover, the trivial bundle is always nef and the tensor product of a nef bundle with an ample bundle is an ample bundle.
Note that \emph{the cotangent bundle of a smooth subvariety $X$ of an abelian variety (of arbitrary characteristic) $A$ is nef}: it is easy to get a surjective map
 $\Omega^1_{A|X}\rightarrow \Omega^1_X$, just consider the exact sequence $$0\longrightarrow \sI_X/\sI_X^2\longrightarrow \Omega^1_{A|X}\longrightarrow \Omega^1_X\longrightarrow 0,$$ where $\sI_X/\sI_X^2$ is the conormal sheaf of $X$ in $A$. Since $$\Omega^1_{A|X}=
 \bigoplus \sO_X$$ is nef then $\Omega_X^1$ is also nef. We define $\Omega_X^p:=\bigwedge ^p\Omega_X^1.$

For subvarieties, the setting is the following.
Let $A$ be a complex abelian variety of dimension $g$ and let $X$
be a complex projective smooth subvariety of $A$ of
dimension $n\geq 2$. Let $M$ be an ample line bundle on $X$ and
$D\in |\omega_X\otimes M|$ be a divisor. Note that the linear
system $|\omega_X\otimes M|$ is non-empty: it is enough to apply
Theorem 5.8 in \cite{FM} to a subvariety of an abelian variety.
Let $$\Sigma_1(D)=\{x\in X \; | \; \mbox{mult}_x D> 1\}$$ be the
singular locus of $D$.

Putting $L:=\omega_X\otimes M$, the sheaf $L^{\otimes s}\otimes \omega_X^{\otimes p}$ satisfies I.T. $0$ for any $s\geq 1$ and any $p\geq 0$: it follows by Kodaira's Vanishing Theorem since $L^{\otimes s}\otimes \omega_X^{\otimes p}= \omega_X\otimes \omega_X^{\otimes {s-1+p}}\otimes M^{\otimes s}$ and the line bundle $\omega_X^{\otimes {s-1+p}}\otimes M^{\otimes s}$ is ample since $\omega_X$ is nef and the tensor product beetwen a nef line bundle and an ample line bundle is still ample. 

\begin{theorem}\label{subvar-cgg}
 The sheaf $\sI_{\Sigma_1(D)}\otimes L^{\otimes s}\otimes \omega_X^{\otimes p}$ is continuously globally generated for any $s\geq 2$ and any $p\geq 1$.

\begin{proof}
Consider the standard exact sequence for the bundle of differential operators of order $\leq 1$ associated to $L$  $$0\longrightarrow \sD_L^0=\sO_X\longrightarrow \sD_L^1\longrightarrow T_X\longrightarrow 0.$$  By tensoring this sequence by $L^{\otimes {s-1}}\otimes \omega_X^{\otimes p}$, with $s\geq 2$ and $p\geq 1$, we get a new one $$0\longrightarrow L^{\otimes {s-1}}\otimes \omega_X^{\otimes p}\longrightarrow \sD_L^1\otimes L^{\otimes {s-1}}\otimes \omega_X^{\otimes p}\longrightarrow \Omega_X^{n-1}\otimes L^{\otimes {s-1}}\otimes \omega_X^{\otimes {p-1}}\longrightarrow 0,$$ where we have used the fact that $\Omega_X^{n-1}=T_X\otimes \omega_X$, see \cite{HA} Exercise II.5.16.

The first term of the sequence satisfies I.T. $0$ and therefore it is continuously globally generated by Proposition \ref{it0 implica cgg}. The surjection $\bigoplus_g \sO_X\rightarrow \Omega_X^1$ of the conormal exact sequence induces a surjection $\bigoplus_{\binom{g}{n-1}} \sO_X\rightarrow \Omega_X^{n-1}$. By tensoring this surjection by $L^{\otimes {s-1}}\otimes \omega_X^{\otimes {p-1}}$, we have that the quotient $\Omega_X^{n-1}\otimes L^{\otimes {s-1}}\otimes \omega_X^{\otimes {p-1}}$ is continuously globally generated by Lemma \ref{cgg-prop} \emph{(i.)}. Now applying Lemma \ref{cgg-prop} \emph{(ii.)} we also get that the middle term of the sequence, $\sD_L^1\otimes L^{\otimes {s-1}}\otimes \omega_X^{\otimes p}$, is continuously globally generated and therefore the quotient $\sI_{\Sigma_1(D)}\otimes L^{\otimes s}\otimes \omega_X^{\otimes p}$ is continuously globally generated.
\end{proof}
\end{theorem}

Proceeding as in the proof of Corollary \ref{last} we easily get the following

\begin{corollary}
 The sheaf $\sI_{\Sigma_1(D)}\otimes L^{\otimes s}\otimes \omega_X^{\otimes p}$ is globally generated for any $s\geq 3$ and any $p\geq 1$.
\end{corollary}

We can also ask for which positive integers $s$ the sheaf $\sI_{\Sigma_1 (D)}$ is cut-out by equations in $|L^{\otimes s}|$. We will use the following vanishing theorem for varieties whose cotangent bundle is nef.

\begin{proposition}\label{vanishing}
 Let $X$ be a complex projective smooth variety of dimension $n$ whose cotanget bundle $\Omega_X^1$ is nef, and let $L$ be an ample
  line bundle on $X$. Then $$H^i(\Omega_X^{p}\otimes \omega_X^{\otimes {p+1}}\otimes L)=0,\quad i>0,\quad p=0,\ldots,n.$$

 \begin{proof}

The cases $p=0,n$ follow directly from the Kodaira's Vanishing Theorem. The idea in general is to apply Demailly's
Vanishing Theorem, see \cite{LA} Theorem 7.3.14. To fix notation, recall briefly the theorem. Given a vector bundle $E$ of
 rank $e$ and a representation $$\rho:GL(e,\mathbb{C})\longrightarrow GL(N,\mathbb{C})$$ of algebraic groups one can
 associate
 to $E$ a bundle $E_\rho$ of rank $N$ by applying $\rho$ to the transition matrices describing $E$.
 The irreducible finite dimensional representations of $GL(e,\mathbb{C})$ are parametrized by non-increasing
 $e$-ples $\lambda=(\lambda_1,\ldots ,\lambda_e)$ where $\lambda_i$ are non negative integers and
 $\lambda_1\geq\ldots\geq\lambda_e\geq 0$. The height of an $e$-ple $h(\lambda)$ is the number of non-zero
  components of $\lambda$. Given $E$ and $\lambda$, we denote by $\Gamma^{\lambda}E$ the bundle associated to the
  representation corresponding to $\lambda$. Note that if $\lambda=(1,\ldots,1,0,\ldots ,0)$ with $m$ repetitions of
   $1$ then $\Gamma ^{\lambda} E=\bigwedge^m E$. Demailly's Vanishing Theorem states that if $E$ is a nef vector bundle
   and $L$ is an ample line bundle then $$H^i(\omega_X\otimes \Gamma^{\lambda}E\otimes (\det E)^{\otimes {h(\lambda)}}\otimes L)=0,
   \quad i>0.$$

Now it is sufficient to apply Demailly's Vanishing Theorem with $E=\Omega_X^1$ and with $\lambda=(\stackrel{p-times}
{\overbrace{1,\ldots,1}},0,\ldots ,0)$, in which case $h(\lambda)=p$.
\end{proof}
\end{proposition}

With the same hypothesis of Theorem \ref{subvar-cgg}, we have 

\begin{theorem}\label{sub-mainth}
The sheaf $\sI_{\Sigma_1(D)}\otimes L^{\otimes s}$ is continuously globally generated for any $s\geq n+2$.

\begin{proof}

By tensoring the standard exact sequence for the bundle of differential operators associated to $L$ of order
$\leq 1$ by $L^{\otimes {s-1}}$, with $s\geq n+2$, we get the following exact sequence $$0\longrightarrow L^{\otimes {s-1}}\longrightarrow \sD_L^1\otimes L^{\otimes {s-1}}\longrightarrow
T_X\otimes L^{\otimes {s-1}}\longrightarrow 0.$$

Let's prove that the bundle $\sD_L^1\otimes L^{\otimes {s-1}}$ satisfies I.T. $0$.\\ The first term of the sequence satisfies I.T. $0$ and by Proposition \ref{vanishing} we get that $$H^i(T_X\otimes L^{\otimes {s-1}}\otimes \alpha)=H^i(\Omega^{n-1}_X\otimes \omega_X^{\otimes n}\otimes \omega_X^{\otimes {s-n-2}}\otimes M^{\otimes {s-1}}\otimes \alpha)=0,$$ $$\forall \;\alpha\in \mbox{Pic}^0(X),\quad i>0,\quad s\geq n+2,$$ therefore also the third term of the sequence satisfies I.T. $0$. Then also the middle term of the sequence satisfies I.T. $0$ and hence it  is continuously globally generated. By Lemma \ref{cgg-prop} \emph{(i.)} the quotient $\sI_{\Sigma_1(D)}\otimes L^{\otimes s}$ is also continuously globally generated.
\end{proof}
\end{theorem}

By Proposition \ref{cgg per cgg=gg} we get the following
\begin{corollary}
 The sheaf $\sI_{\Sigma_1(D)}\otimes L^{\otimes s}$ is globally generated for any $s\geq n+3$.
\end{corollary}


\section*{Acknowledgemets}
It is a pleasure to thank Giuseppe Pareschi and Mihnea Popa for many valuable discussions and for their lectures held during the summer school PRAGMATIC 2007 in Catania (Italy). Furthermore we want to thank the Department of Mathematics of University of Catania for the nice stay, where this project was started.

\small{\texttt{UNIVERSITA' DEGLI STUDI DI ROMA, TOR VERGATA,\\
VIA DELLA RICERCA SCIENTIFICA, 1 - 00133 ROMA}}\\
\textit{email:}\texttt{lombardi@mat.uniroma2.it}\\

\small{\texttt{POLITECNICO DI TORINO,\\
CORSO DUCA DEGLI ABBRUZZI, 24 - 10129 TORINO}}\\
\textit{email:}\texttt{malaspina@calvino.polito.it}\\

\end{document}